\documentclass{amsart}
 \usepackage{amsfonts}
 \begin{document}
\newtheorem{Theorem}{Theorem}[section]
\newtheorem{Remark}{Remark}[section]
\newtheorem{Lemma}{Lemma}[section]
\newtheorem{Proposition}{Proposition}[section]
\newtheorem{Corollary}{Corollary}[section]
\def\qedbox{\hbox{$\rlap{$\sqcap$}\sqcup$}}
\def\qed{\nobreak\hfill\penalty250 \hbox{}\nobreak\hfill\qedbox}
 \def\e{\varepsilon}
 \def\diag{\mbox{diag\,}}
 \def\dim{\mbox{dim\,}}
 \def\arg{\mbox{arg\,}}
 \def\nnull{\mbox{null\,}}
 \def\tr{\mbox{tr\,}}
 \def\xt{\tilde{X}}
  \def\yt{\tilde{Y}}
  \def\zt{\tilde{Z}}
  \def\wt{\tilde{W}}
  \def\a{\alpha}
  \def\b{\beta}
  \def\oo{\omega}
 \def\e{\epsilon}
 \def\ve{\varepsilon}
 \def\RR{\mathbb{R}}
 \def\CC{\mathbb{C}}
 \def\HH{\mathbb{H}}
 \def\FF{\mathbb{F}}
 \def\NN{\mathcal{N}}
 \def\ZZ{\mathbb{Z}}
 \def\II{\mathbb{1}}
 \def\g{\mathfrak g}
\def\ggg{\mathfrak g}
\def\rrr{\mathfrak r}
\def\sss{\mathfrak s}

\def\ej{e^1}
\def\ed{e^2}
\def\et{e^3}
\def\ec{e^4}

\def\ejj{e_1}
\def\edd{e_2}
\def\ett{e_3}
\def\ecc{e_4}

\title{ Four-dimensional Lie algebras with a para-hypercomplex structure}
\author{Novica Bla{\v z}i\' c and  Srdjan Vukmirovi\' c}
\address{Faculty of mathematics, University of Belgrade, Studenski
trg 16, p.p. 550, 11 000 Belgrade, Yugoslavia}
\email{blazicn@matf.bg.ac.yu,  vsrdjan@matf.bg.ac.yu}
\thanks{Research partially  supported by MNTS, project $\#$1854.}
\date{\today\ phcs19}
\begin{abstract}
The  main goal  is to classify $4$-dimensional real  Lie algebras
$\g$ which admit a para-hypercomplex structure. This is a step
toward the classification  of Lie groups admitting the
corresponding left-invariant structure and therefore possessing  a
neutral, left-invariant, anti-self-dual metric. Our study is
related to the work of Barberis who classified real,
$4$-dimensional simply-connected Lie groups which admit an
invariant hypercomplex structure.
\end{abstract}
\maketitle
\section{Introduction}
Our work is motivated  by the work of Barberis~\cite{B} where
invariant hypercomplex structures  on $4$-dimensional real   Lie
groups are classified (see Section \ref{sec:prels} for
definitions). In that case the corresponding hermitian metric is
positive definite and unique up to a positive constant. Our main
goal is to classify $4$-dimensional real  Lie algebras $\g$ which
admit para-hypercomplex structures. This is a step toward the
classification of the corresponding left invariant structures on
Lie groups. In this case  the corresponding hermitian
pseudo-Riemannian metric   determined  by the para-hypercomplex
structure is also  unique up to a constant, but has to be of
signature $(2,2)$. This metric is anti-self-dual (see
~\cite{Sal}).

In the paper \cite{Sal1}  Andrada and Salamon have shown that any
para-hypercomplex structure on a real Lie algebra $\g$ rise to a
hypercomplex structure on its complexification $\g ^\CC$
(considered as a real Lie algebra). They referred to
para-hypercomplex structure as complex product structure.

Let us remark that Snow~\cite{S} and Ovando~\cite{Ov}  classified
the invariant complex structures on $4$-dimensional, solvable,
simply-connected real  Lie groups where the  dimension of
commutators is less than three and equal three, respectively.
Since every para-hypercomplex manifold is also complex, the Lie
algebras from our classification also appear in their lists.

Let us state the main theorem (proved in Subsection
\ref{ssec:exist}).

\begin{Theorem}
\label{th:main} Up to an isomorphism the only $4$-dimensional Lie
algebras $\ggg$ admitting an integrable
 para-hypercomplex structure  are listed below.
\begin{itemize}
\item[(PHC1)] $\g$ is abelian,
\item[(PHC2)] $[X, Y]=W, \enskip [Y, W]=-X, \enskip [W, X]=Y,$
\item[(PHC3)] $[X,Y] = Y, \enskip [X, W]
= W,$
\item[(PHC4)] $[X, Y] = Z,$
\item[(PHC5)] $[X, Y] = X,$
\item[(PHC6)]
$[X, Y] = Z, \enskip [X, W] =  X + aY + bZ, \enskip [W, Y] = Y,$
\item[(PHC7)] $[X,Z]=X, \enskip [X,W]=Y, \enskip [Y, Z]=Y, \enskip
[Y,W]=aX+bY,\enskip
           a,b\in \RR,$
\item[(PHC8)] $[X,Z]=X,\enskip [Y,W]=Y,$
\item[(PHC9)]{ $ [Z, W] = Z, \enskip  [Y, W] = Y, \enskip [X, W] = cX + aY + bZ,\, c\neq 0, \enskip  a, b \in \RR ,$}
\item[(PHC10)]{$[Y, X] =  Z, \enskip [W, Z] = cZ,\enskip  [W, X] = \frac{1}{2}X + aY + bZ, [W, Y] = (c- \frac{1}{2}   ) Y, c \neq 0.$ }
\end{itemize}
In the previous list  the additive basis of algebra $\ggg$ is $(X,
Y, Z, W)$, and only the non-zero commutators are given.
\end{Theorem}

In the proof we study separately the cassis defined in the terms
of  metrics defined on the  derived algebra $\ggg '$ by means of
the para-hypercomplex structure.

Here is a brief outline of the paper. In Section \ref{sec:prels}
we first give necessary definitions and prove some basic
properties of para-hypercomplex structures and  a number of lemmas
which we use in the sequel. In Section \ref{sec:algebras} we
step-by-step prove Theorem \ref{th:main}. First, in Subsection
\ref{ssec:withcenter} we classify $4$-dimensional Lie algebras
with a non-trivial center and admitting a para-hypercomplex
structure. Further on we suppose that algebra $\ggg$ has a trivial
center. In Subsection \ref{ssec:dim2} and \ref{ssec:dim3} we
classify solvable $4$-dimensional  Lie algebras $\ggg$ admitting a
para-hypercomplex structure (Theorems \ref{th:dim1}, \ref{th:dim2}
and \ref{th:g3} depending on  the dimension of the commutator
subalgebra $\ggg ' = [\ggg, \ggg]$). In Subsection
\ref{ssec:exist} we prove the Theorem \ref{th:main} using previous
classifications and find  particular examples of para-hypercomplex
structures on algebras PHC1-PHC10. Finally, in Section
\ref{sec:comp} we compare our results with the results of
Barberis~\cite{B}.

\section{Preliminaries}
\label{sec:prels}

Let $V$ be a real vector space. A {\it complex structure} on $V$
is an endomorphism $J_1$ of $V$ satisfying the condition
$$
J_1^2 = -1.
$$
Existence of a complex structure implies that $V$ has to be of an
even dimension. A {\it product structure} on $V$ is an
endomorphism $J_2$ of $V$ satisfying the conditions
$$
J_2^2 = 1, \quad J_2\neq \pm 1.
$$
 A {\it para-hypercomplex  structure} on  $V$ is a pair $(J_1, J_2)$ of anti-commuting
complex structure $J_1$ and product
 structure $J_2$, i.e.
 satisfying the relations
\begin{equation}
\label{eq:comm_rel} J_1^2 = -1, \enskip J_2^2 = 1, \enskip J_1J_2
= -J_2 J_1.
\end{equation}
If both structures $J_1$ and $J_2$ are complex then the pair
$(J_1, J_2)$ is called a {\it hypercomplex structure} on $V.$ In
the sequel we concentrate on the case of para-hypercomplex
structure.

It is customary to denote $J_3 = J_1J_2.$ Note that the structure
$J_3$ is a product structure. The Lie subalgebra of $End(V)$
spanned by $J_1, J_2$ and $J_3$  is isomorphic to $sl_2(\RR).$ Any
$x= (x_1, x_2, x_3)\in \RR ^3$ defines a  structure by the formula
\begin{equation*}
\label{eq:Jx} J_x :=  x_1 J_1 + x_2J_2 + x_3J_3 .
\end{equation*}
Denote by
$$
 \langle x  , y \rangle = x_1y_1-x_2y_2-x_3y_3,
$$
$x = (x_1, x_2, x_3), \enskip  y = (y_1, y_2, y_3),$  the  inner
product
 in $\RR ^3 = \RR ^{1,2}$ and by
$$
x\times y = (x_2y_3 - x_3y_2,  x_3y_1 - x_1y_3, x_1y_2 - x_2y_1)
$$
the usual cross product.
 The structure $J_x$ is  a complex structure provided that
$$
\langle x,x \rangle  = x_1^2- x_2^2 - x_3^2 = 1
$$
and a product structure provided that
$$
\langle x,x \rangle  = x_1^2- x_2^2 - x_3^2 = -1.
$$
Hence, a para-hypercomplex structure $(J_1, J_2)$ defines a
$2$-sheeted  hyperboloid ${\ss}^-$ of complex structures and a
$1$-sheeted hyperboloid ${\ss}^+$ of product structures.
\begin{Proposition}
If $(J_1, J_2)$ is a para-hypercomplex structure on a vector space
$V,$ then:
\begin{itemize}
\item[i)]{$ J_x J_y = -\langle x, y \rangle 1 + J_{x\times y}.$}
\item[ii)]{The pair $(J_x, J_y)\in {\ss} ^- \times {\ss} ^+$ is a
para-hypercomplex structure if and only if $x \perp y.$}
\end{itemize}
\end{Proposition}
\noindent {\bf Proof:}  From  the relations
$$
J_1J_2 = J_3 = -J_2J_1, \quad J_1J_3 = -J_2 = -J_3J_1, \quad
J_2J_3 = -J_1 = -J_3J_2
$$
 the statement i) follows by a direct calculation.

 Since $J_x$ is a complex structure and $J_y$ is a product
structure, the pair $(J_x, J_y)$ is a para-hypercomplex structure
if and only if $J_x$ and $J_y$ anti-commute. Using the relation i)
and the anti-commutativity of the cross product we have
$$
0 = J_xJ_y + J_yJ_x = -2\langle x, y\rangle 1.
$$
Hence, the statement ii) is proved.\qed

\noindent The para-hypercomplex structures $(J_1, J_2)$ and $(J_x,
J_y)$  are called {\it compatible}. An {\it almost
para-hypercomplex structure} on a manifold $M$ is a pair $(J_1,
J_2)$ of sections of $End(TM)$ satisfying the relations
(\ref{eq:comm_rel}). It is a para-hypercomplex structure if both
structures are {\it integrable}, that is, if the corresponding
Nijenhuis tensors
\begin{equation}
\label{eq:nijenhuis}
 \NN_\alpha (X, Y) = [J_\alpha X, J_\alpha Y]-
J_\alpha [X, J_\alpha Y] - J_\alpha [J_\alpha X, Y] \pm [X, Y],
\end{equation}
$\alpha = 1,2,$ vanish on all  vector fields $X, Y$. In this
formula sign $-$ occurs in the case of a complex structure and
sign $+$ occurs in the case of a product structure.

If $M=G$ is a Lie group we additionally assume that the
para-hypercomplex structure is left invariant. This allows us to
also  describe a para-hypercomplex structure on its Lie algebra
$\g$.
 Hence, a para-hypercomplex structure  $(J_1, J_2)$ on
$\g$ satisfies  both relations  (\ref{eq:comm_rel}) and
(\ref{eq:nijenhuis}).

\begin{Proposition}
Let $(J_1, J_2)$ be an integrable  para-hypercomplex structure on
a Lie algebra $\g.$
\begin{itemize}
\item[i)]{The product structure $J_3 = J_1J_2$ is integrable.}
\item[ii)]{Any compatible para-hypercomplex structure $(J_x, J_y)$ is
integrable.}
\end{itemize}
\end{Proposition}
\noindent {\bf Proof:}  The statement i) follows from the relation
\begin{eqnarray*}
2 \NN _3(X, Y) &=& \NN _1(J_2X, J_2Y) + \NN _2(J_1X, J_1Y) -
J_1\NN _2(J_1X, Y) - J_1\NN _2(X, J_1 Y) +\\
&+& \NN _2 (X, Y) - J_2 \NN _1 (J_2 X, Y) - J_2 \NN _1(X, J_2Y ) -
\NN _1 (X, Y)
\end{eqnarray*}
where $\NN _3$ is the Nijenhuis tensor of the product structure
$J_3.$

To prove ii) denote by $\NN _x$ the Nijenhuis tensor corresponding
to the structure $J_x, \, x = (x_1, x_2, x_3).$ One can check that
\begin{eqnarray*}
\NN _x  &=& x_1^2\NN _1 + x_2^2\NN _2 + x_3^2 \NN _3 +
x_1x_2(J_3\NN _1 + J_3\NN _2
+ J_3\NN _3J_1) +\\
&+& x_2x_3(J_1\NN _2- J_1\NN _3-J_1\NN _1J_2) + x_1x_3(-J_2\NN
_1-J_2\NN _3+J_2\NN _2J_3)
\end{eqnarray*}
holds, where we have used the notation, for instance
$$
J_2\NN _2J_3 (X, Y) = J_2\NN _2(J_3X, J_3Y).
$$
Now, statement ii) follows using statement i).\qed

Let $g$
be an inner product on the vector space $V$. A para-hypercomplex
structure $(J_1, J_2)$ on $V$  is called hermitian with respect to
$g$ if
\begin{equation}
\label{eq:hermitian} g( J_\alpha X, Y) = -g( X, J_\alpha Y ) ,
\quad X, Y \in V
\end{equation}
holds, i.e. if both structures $J_1$ and $J_2$ are hermitian. It
is easy to prove that a hermitian complex structure is an isometry
and a hermitian product structure is an anti-isometry, i.e.
\begin{equation*}
\label{eq:anti_isometry} g( J_1 X, J_1Y) = g( X, Y), \quad g( J_2
X, J_2Y ) = -g( X, Y).
\end{equation*}
Existence of an anti-isometry  implies that the inner product $g$
must be  of  neutral, $(n,n)$ signature.
\begin{Proposition}
\label{cor:metric}
 Let $(J_1, J_2)$ be a para-hypercomplex
structure hermitian with respect to the scalar product $g$ on the
vector space $V.$
\begin{itemize}
\item[i)]{The product structure $J_3 = J_1J_2$ is hermitian.}
\item[ii)]{Any compatible para-hypercomplex structure $(J_x, J_y)$ is
hermitian.}
\end{itemize}
\end{Proposition}
\noindent {\bf Proof:} i) If  $J_1$ and $J_2$ are hermitian then
$J_3$ is hermitian since we have
$$
\langle J_3 X, Y\rangle = \langle J_1J_2 X, Y\rangle = -\langle
J_2 X, J_1Y\rangle = \langle  X, J_2 J_1Y\rangle  =-\langle  X,
J_3Y\rangle.
$$
ii) Since the condition of any $J_x$ to be hermitian is linear
with respect to $x$, the statement ii) follows from the statement
i).\qed

Now, we prove some lemmas which will be useful in the sequel.

\begin{Lemma}
\label{le:unique_product} If $(J_1, J_2)$ is a para-hypercomplex
structure on a
 real $4$-dimensional vector space $V$ then:
\begin{itemize}
\item[i)]
 There is an
inner product $g$ on $V$, unique up to a non-zero constant, such
that the structure $(J_1, J_2)$ is hermitian with respect to $g$.
\item[ii)] Any compatible para-hypercomplex structure $(J_x, J_y)$
determines the same inner product $g$ on $V.$
\end{itemize}
\end{Lemma}
\noindent {\bf Proof:} First, we prove the existence of such an
inner product. If $(\cdot, \cdot)$ is an arbitrary inner product
on $V,$ then the
 inner product
\begin{equation}
g( X, Y) := (X, Y)  + (J_1X, J_1Y) - (J_2X, J_2Y)-(J_3X, J_3Y)
\end{equation}
satisfies the properties (\ref{eq:hermitian}).

To see the uniqueness let $g'( \cdot, \cdot )$ be another inner
product on $V$ satisfying (\ref{eq:hermitian}). As remarked before
 both products are of  signature $(2,2)$. There exists a vector
$X$ which is not null with respect to the both inner products, for
instance
\begin{equation*}
g( X, X)  = 1, \quad g'( X , X ) = \lambda \neq 0.
\end{equation*}

The relations  (\ref{eq:comm_rel}) and (\ref{eq:hermitian}) imply
that the vectors $X,$ $J_1 X,$ $J_2X,$ $J_3X$ are mutually
orthogonal with respect to both inner products. Moreover,
\begin{eqnarray*}
g( X, X) = g( J_1X,J_1X) = &1&
  = -g(J_2X, J_2X) = -g( J_3X, J_3X)\\
g'( X, X) = g'( J_1X, J_1X)  = &\lambda& = -g'( J_2X, J_2X) = -g'(
J_3X, J_3X).
\end{eqnarray*}
Hence, $g(\cdot , \cdot ) = \lambda g'( \cdot,\cdot ),
  \enskip \lambda \neq 0.$

ii) According to Proposition \ref{cor:metric} the structure $(J_x,
J_y)$ is hermitian with respect to $g$. The statement follows from
the uniqueness of $g$ (up to a non-zero scalar).\qed
\begin{Remark}
\label{re:unique_product}
 In the light of Lemma
\ref{le:unique_product} we see that the notion of null vector $N$
(such that $g( N, N)= 0$) depends only on the hermitian structure
$(J_1, J_2)$ and not on a particular inner product.
\end{Remark}
From the proof of  Lemma \ref{le:unique_product} we also obtain
the following.
\begin{Lemma}
\label{le:J_basis} If $(J_1, J_2)$ is a is a para-hypercomplex
structure on a  real $4$-dimensional vector space $V$  then
$$
(X, J_1X, J_2X, J_3X)\enskip  \mbox{is a basis of } V \quad
\Leftrightarrow\quad X \mbox{ is not null.}
$$
\end{Lemma}
\begin{Lemma}
\label{le:int_cond}
 If $J_\alpha$ is an endomorphism of a $4$-dimensional Lie
algebra $\g$ such that $J_\alpha ^2 = \pm 1$ and $(X, J_\alpha X,
Y, J_\alpha Y)$ is a basis of $\g$ then the corresponding
Nijenhuis tensor  $N_\alpha$ vanishes if and only if $N_\alpha (X,
Y) = 0.$
\end{Lemma}
\noindent {\bf Proof:} One can easily show that $\NN _\alpha
(J_\alpha X, Y) = - J_\alpha \NN _\alpha (X, Y).$ The lemma
follows from the fact that $\NN _\alpha$ is antisymmetric and
bilinear .\qed
\begin{Lemma}
\label{le:2d_plane} Let $(J_1, J_2)$ be  a para-hypercomplex
structure on a  real $4$-dimensional vector space $V$ and let
$W\subset V$ be a $2$-dimensional subspace. Then there exists a
compatible para-hypercomplex structure $(J_1', J_2')$ such that:
\begin{itemize}
\item[i)]{If $W$ is  definite (contains no null directions)
 then $J_1'W = W.$}
\item[ii)]{If $W$ is Lorentz (contains exactly two null directions) then $J_2'W
= W.$}
\item[iii)]{If $W$ is totally null (every  vector in $W$ is a null vector)
then either\\
(a) $\quad J_2' |_W = 1, \quad  V = W \oplus J_1'W, \quad
\mbox{or}$\\
 (b) \, there exists a non-null vector $X$ such that
$$
\quad \quad \quad \quad W = \RR \langle J_1'X + J_2'X, X - J_3'X
\rangle , \quad
 J(W)= W \enskip \mbox{ for all } \enskip  J\in{\ss}^\pm .
$$ }
\item[iv)]{If the induced metric on $W$  is of rank $1$ ($W$ contains
exactly one null direction $N$) then  $N = J_1'X -J_2'X$  for any
given vector $X \in  W, \, |X|^\neq 0.$ }
\end{itemize}
\end{Lemma}
\noindent {\bf Proof of  i) and ii):} Let $(X, Y)$ be a
pseudo-orthonormal basis of $W$ ($|X|^2 = -|Y|^2 = 1$ and $\langle
X, Y\rangle = 0$ with respect to the induced inner product on
$W$). Then, according to  Lemma \ref{le:J_basis} vectors $X,$
$J_1X$, $J_2X$ and $J_3X$ form a pseudo-orthonormal basis of $V$
and  we have $Y =x_1J_1X + x_2J_2X + x_3J_3X$ with $x_1^2 - x_2^2
- x_3^2 = \pm 1 ,$ where $-$ occurs if $W$ is Lorentz and $+$ if
$W$ is positive or negative definite. The structure
\begin{equation*}
J_x = x_1J_1 + x_2J_2 + x_3J_3
\end{equation*}
preserves $W$. It is a product structure if $W$ is Lorentz (and we
set $J_2' =J_x$) and a complex structure if $W$ is definite (and
we set $J_1' =J_x$). The second structure can be chosen such that
$(J_1', J_2')$ forms a compatible para-hypercomplex structure.
Note that there cannot exists a product structure preserving a
definite $W$ since a product structure is an anti-isometry.
Similarly, a complex structure preserving a Lorentz $W$ cannot
exist.

\noindent {\bf Proof of  iii)} Let $N_1\in W$ be a null vector.
There exists a non-null vector $X\in V$ perpendicular to $N_1$.
Hence
$$
N_1 = \alpha  J_1X + \beta J_2X + \gamma J_3X
$$
and
$$
\alpha ^2 - \beta ^2 - \gamma ^2 =0
$$
so $\alpha  \neq 0$ and we may assume that $\alpha =1$. Then
$J_2'= \beta  J_2+ \gamma J_3$ is a product structure, the
structure  $(J_1', J_2'), \enskip J_1' = J_1$ is a compatible
para-hypercomplex structure and we have
$$
N_1 = J_1'X+ J_2'X.
$$
Any null vector $aX + bJ_1'X + cJ_2'X + dJ_3'X$ which is
orthogonal to the vector $N_1$ is of the form
$$
N^\pm = aX + bJ_1'X + bJ_2'X \pm aJ_3'X.
$$
Notice that the vector $N_1$ is also of the form $N^\pm$ and that
there exist exactly two null planes $W^\pm$ containing the vector
$N_1.$ They can be written in the form
$$
W^\pm = \RR \langle N_1 , N_2^\pm = X \pm J_3'X \rangle .
$$
The plane $W^-$ is  the  $+1$-eigenspace of the product structure
$J_3'$ and  the vectors $N_1,$ $N_2^-,$ $J_1'N_1,$ $J_1'N_2^-$ are
independent, so $V = W^- \oplus J_1'W^-$ and iii)a holds.

In the case of the plane $W^+$ one easily checks that
 $J_1'W^+ = W^+ = J_2' W^+$ and hence statement iii)b follows.

\noindent {\bf Proof of  iv)}  The proof is similar to the first
part of the previous  proof (with $N_1=N$).  \qed

\begin{Lemma}
\label{le:3d_plane}
Let $(J_1, J_2)$ be  a para-hypercomplex
structure on a  real $4$-dimensional vector space $V$ and let
$W\subset V$ be a $3$-dimensional subspace such that the induced
metric is degenerate.
For $N\in W^\perp$
and  $X\in W$, $|X|^2\neq 0$,
there exists a  compatible  para-hypercomplex structure
 $(J_1', J_2')$ on $V$  such that
$N=J_1'X-J_2'X$ and the arbitrary null vector in $W$ belongs to
the union of two-dimensional planes $ \pi_1=\RR\langle
N,J_1'N\rangle$  and $\pi_-= \{V\mid J_3'V=-V\}, $ i.e.
$$
\nnull(W)=\{U\in W\mid |U|^2=0\}= \pi_1\cup\pi_- =\RR\langle
N,J_1'N\rangle \cup \{V\mid J_3'V=-V\}.
$$
\end{Lemma}
\noindent {\bf Proof:} Since we have $|N|^2 =0, \, |X|^2 \neq 0,
\, \langle  N, X\rangle  = 0$ the existence of a compatible
structure $(J_1', J_2')$ such that $N = J_1'X - J_2'X$ follows
from the Lemma \ref{le:2d_plane}  iv). Moreover,
 $\{N,J_1'N,X\}$ is a basis of $W$ and
$\{N,J_1'N,X, J_1'X\}$ is a basis of $V$. Thus, for $U\in \nnull
(W)$ of the form  $U=\alpha N+\beta J_1'N +\gamma X$ we get
$$
0 =|U|^2=\gamma(\gamma-2\beta)|X|^2.
$$
The case  $\gamma = 0$ gives the plane $\pi _1 = \RR \langle N,
J_1N \rangle.$ For $\gamma = 2\beta$ one can check that $J_3(U) =
-U$, so $U$ belongs to the $-1$ eigenspace of $J_3'.$ \qed

\section{Lie algebras admitting a para-hypercomplex structure}
\label{sec:algebras}
\subsection{Case when $\g$  has a non-trivial center}
\label{ssec:withcenter}
 In the following theorem the additive
basis of the Lie algebra $\g$ is either $(X, Y, Z, W)$ or $(X,Y,
N_1, N_2).$ The vectors $N_\alpha$ are null vectors.
\begin{Theorem}
\label{th:withcenter}
 A $4$-dimensional Lie algebra $\g$ admitting a para-hypercomplex structure and
with a non-trivial center $Z(\g )$ is one of algebras PHC1-PHC6.
\end{Theorem}
As a consequence of Levi decomposition theorem and the
classification of real semisimple Lie algebras the only
non-solvable Lie algebras which are $4$-dimensional are $\RR
\oplus so(3)$ and $\RR \oplus sl_2(\RR )$. Since they both have a
non-trivial center, as a consequence of  Theorem
\ref{th:withcenter} we have the following corollary.
\begin{Corollary}
The only non-solvable, real $4$-dimensional Lie algebra admitting
a para-hypercomplex structure  is $\RR \oplus sl_2(\RR ).$
\end{Corollary}
\noindent {\bf Proof of Theorem \ref{th:withcenter}:} In order to
prove that these are the only Lie algebras with non-trivial center
which  admit a para-hypercomplex structure we consider two cases.

{\bf Case 1:} there exists a non-null central element $Z.$
\noindent Let $(J_1, J_2)$ be a para-hypercomplex structure on
$\g$ and denote
\begin{equation*}
X = J_1Z,\enskip Y = J_2Z, \enskip W = J_3Z.
\end{equation*}
Then
\begin{equation}
\label{eq:[X,Y]}
 [X, Y] = a Z + b X + cY + dW.
\end{equation}
According to  Lemma \ref{le:int_cond} integrability of $J_1$ is
equivalent to
\begin{equation}
\label{eq:int1}
 0 = \NN_1(Z, Y) = [X, W] - J_1[X, Y].
\end{equation}
Similarly, the integrability of $J_2$ is equivalent to
\begin{equation}
\label{eq:int2}
 0 = \NN_2(X, Z) = [Y, W] - J_2[X, Y].
\end{equation}
From the relations (\ref{eq:[X,Y]}), (\ref{eq:int1}) and
(\ref{eq:int2}) we get
\begin{equation*}
[X, W] = -bZ+ aX  - dY+ cW,\quad [Y, W]
 =  cZ - dX + a Y - b W.
\end{equation*}
The Jacobi identity  is equivalent to
\begin{eqnarray*}
0&=& [[X, Y],W] + [[Y, W],X]+ [[W,X],Y] =\\
&=& 2(-a^2 - b^2 + c^2)Z -2cd X - 2dbY - 2ad W.
\end{eqnarray*}
If $a=b=c=d=0$ then the algebra $\g$ is abelian, i.e. PHC1.
\noindent If $a=b=c=0$ and $d\neq 0$ then after scaling $\g \cong
R\oplus sl_2(\RR )$, i.e. PHC2.

 \noindent If $d = 0$ and $0 \neq
c^2 = a^2 + b^2$ then the derived algebra $\g ' = [\g , \g]$ of
$\g$ is $2$-dimensional since
\begin{equation*}
c[Y,W] = a[X,Y] + b[W,X]
\end{equation*}
It is generated by the vectors $W_1 = [X, Y],  \enskip Y_1 = [W,
X].$ The vectors $Z,$ $X_1 = \frac{1}{c}X,$ $Y_1$ and $W_1$ are
linearly independent and we get algebra PHC3.

{\bf Case 2:} all central vectors are null vectors. Denote one of
them by $N$.  According to  Lemma \ref{le:2d_plane} iv),  we can
assume that $N = J_1X-J_2X$ for a non-null vector $X\in \g '$.
Then the vectors $N, J_1N, X$ and $J_1X$ form a basis of $\g$ and
the structure $J_2$ expressed in the terms of that basis reads
\begin{equation}
J_2X = J_1X-N, \enskip J_2J_1N = N, \enskip J_2J_1 X = J_1N +
X,\enskip J_2N = J_1N.
\end{equation}
The integrability of the structure $J_1$ gives the following
conditions
\begin{equation}
\label{eq:int11}
 0 = \NN_1(X, N) = [J_1X, J_1N] - J_1[X, J_1N].
\end{equation}
Since the vectors $N, J_2N, X$ and $J_2X$ form a basis of $\g$,
the
integrability of the product structure $J_2$ is equivalent to
\begin{equation}
\label{eq:int21}
 0 = \NN_2(X, N) = [J_1X, J_1N] - J_2[X, J_1N].
\end{equation}

The vector $[X, J_1N]$ is of the  form $[X, J_1N] = aN + bJ_1N +
cX + dJ_1X.$ Using the relations (\ref{eq:int11}) and
(\ref{eq:int21}) we get that
\begin{equation}
[X, J_1N]  = aN + bJ_1N + 2b X, \quad [J_1X, J_1N ]= -bN + aJ_1N +
2bJ_1X.
\end{equation}
If we write $[X, J_1X] = \alpha N + \beta J_1N + \gamma x + \delta
J_1X$ and impose the Jacobi identity on the vectors $J_1N, X$ and
$J_1X$ we get the following system of equations:
\begin{eqnarray*}
-4\alpha b - b^2 - \delta b + \gamma a - a^2& =& 0,\\
- 4 b \beta + a \delta + b \gamma &=& 0,\\
b(a+ \gamma ) &=& 0,\\
b(b- \delta) &=& 0.
\end{eqnarray*}
The system has three classes of solutions.

 \noindent {\bf 1) $a =
0 = b.$} In this case the only non-zero commutator is
\begin{equation*}
[X, J_1X] = \alpha N + \beta J_1N + \gamma X + \delta J_1X.
\end{equation*}
If $\gamma = 0 = \delta$, the change of the basis $Y = J_1X,$ $N_1
= \alpha N + \beta J_1N,$ $N_2 \in \RR \langle N, J_1N\rangle$
 gives the relations PHC4.
\noindent If $\delta \neq 0$ then the change $Y =
\frac{1}{\delta}[X, J_1X],$ $N_1 = N, N_2 = J_1N$ gives the
relations PHC5. The case $\delta = 0, \gamma \neq 0$ similarly
reduces to the relations PHC5.

\noindent {\bf 2) $ b = \delta \neq 0,a = -\gamma.$} This case
reduces to the relations PHC3.

\noindent {\bf 3) $a = \gamma \neq 0.$} This immediately gives the
commutator relations PHC6. \qed

\subsection{Case of solvable Lie algebra $\g$  and $\dim \g'\leq 2$.}
\label{ssec:dim2}

\begin{Theorem}
\label{th:dim1}
 Let $\g $ be a $4$-dimensional real Lie algebra
admitting a para-hy\-per\-com\-plex structure and $\dim \g '=1.$
Then  $\g$ is one of the algebras PHC1, PHC2 from Theorem
\ref{th:withcenter}.
\end{Theorem}
\noindent {\bf Proof:} If $\g$ has a non-trivial center $\xi$ then
from  Theorem \ref{th:withcenter} we get the algebras PHC1 and
PHC2. Now, as in~\cite{B}, Proposition 3.2, let $\xi = \{0\}$ and
let $X$ be a non-zero element of $\g'.$ There exists $Y$ such that
$[Y,X] = X.$ Then $\g$ decomposes as
$$
\g = ker (ad_X) \cap ker (ad_Y) \oplus \RR X \oplus \RR Y.
$$
From the Jacobi identity we get that $\xi = ker (ad_X) \cap ker
(ad_Y)$, a contradiction. Hence solvable $\g$ without center and
with $\dim \g ' = 1$ does not exist (this does not depend on the
existence of para-hypercomplex structure). \qed

\begin{Theorem}
\label{th:dim2}
 \label{th:withoutcenter}
 Let $\g$ be a $4$-dimensional solvable
Lie algebra admitting a para-hypercomplex structure and with $\dim
\g '=2.$ If $\g$ has a non-trivial center then it is algebra PHC2.
If $\g$ has a trivial center then $\g$ is one of algebras
PHC7-PHC9.
\end{Theorem}
\begin{Remark}
\label{re:dim2}
 Using the notation introduced by Snow~\cite{S},
these Lie algebras are S11, S8 and S10 respectively. The class S11
contains
 as a special case the Lie algebra $\frak{aff}(\CC)$
which is the unique solvable Lie algebra with $2$-dimensional
derived algebra which admits hypercomplex structure \cite{B}.
\end{Remark}
\noindent {\bf Proof:} Suppose that the center of $\g$ is trivial
and that $(J_1, J_2)$ is a para-hypercomplex structure on $\g.$
According to Lemma \ref{le:unique_product} and Remark
\ref{re:unique_product} the structure $(J_1, J_2)$ determines the
inner product on $\g = V$ and the notion of a null vector.
 As in Lemma~\ref{le:2d_plane} we have to consider
 the cases concerning the rank and the signature of the
induced inner product  on $\g'=W$.

{\bf  Case i):} Induced metric on $\g'$ is definite. \noindent
  Because of  Lemma~\ref{le:2d_plane}  i) we may assume  that
  $\g'$ is invariant with respect to the complex structure $J_1$,
$J_1\g'=\g'$,
and $\g=\g'\oplus J_2 \g'$.
  Let  $\{X, J_1X=Y\}$ be  a basis of $\g'$ and
  $\{X, Y,J_2X,J_2Y\}$ be a basis of $\g$.
  The Lie algebra $\g'$ is abelian since $\g$ is solvable and by the
integrability
of the product  structure $J_2$ we have  $\NN_2(X,J_1X)=0$ and
\begin{equation}
\label{eq:d2d2}
[J_2X, J_2Y]=0, \quad [J_2X,Y]=[J_2Y,X].
\end{equation}
Because of the integrability of the complex structure $J_1$,
$\NN_1(X,J_2X)=0$ and
\begin{equation}
\label{eq:d2d1}
[X, J_2X]=-[Y, J_2Y].
\end{equation}
For arbitrary vectors $V$ and  $W$ in $\g$,
$$
 [V,W]=\a(V,W)X+\b(V,W)Y,
$$
where $\a$ and $\b$ are skew-symmetric bilinear forms on $\g$.
From the Jacobi identity we have
$$
   \a(X,J_2X)=\b(X,J_2Y), \quad \a(J_2Y,X)=\b(X,J_2X)
$$
and the bracket in $\g$ is determined by $c=\a(X,J_2X)$ and
$d=\b(X,J_2X)$ as follows:
$$
[X, J_2X]=-[Y, J_2Y]=cX+dY,  \quad
[X, J_2Y]=[Y, J_2X]=-dX+cY.
$$
Since $\dim \g'=2$,  $c^2+d^2\neq 0$ and we may choose
\begin{eqnarray*}
\xt = (c^2+d^2)^{-1}(cX+dY),{\phantom{J_2J_2}}  \qquad &~&
\yt=(c^2+d^2)^{-
1}(-dX+cY),\\
\zt = (c^2+d^2)^{-1}(cJ_2X-dJ_2Y), \qquad  &~& \wt=
(c^2+d^2)^{-1}(dJ_2X+cJ_2Y),
\end{eqnarray*}
and hence
\begin{eqnarray*}
~[\xt,\zt]=\xt, \qquad &~& [\xt,\wt]=\yt, \\
~[\yt,\zt]=\yt, \qquad &~& [\yt,\wt]=-\xt,
\end{eqnarray*}
so we get the algebra PHC7 for $a=-1, b=0$. Note that  $\g\equiv
\frak{aff}(\CC)$.

{\bf Case ii):} Induced metric on $\g'$ is indefinite, of Lorentz
type $(-+)$. \noindent Because of  Lemma~\ref{le:2d_plane}  ii) we
may assume that $\g'$ is invariant with respect to the product
structure $J_2$, $J_2\g'=\g'$, and $\g=\g'\oplus J_1 \g'$. Let
$\{X, J_2X=Y\}$ be a basis of $\g'$ and  $\{X, Y,J_1X,J_1Y\}$ be a
basis of $\g$. By the integrability of the complex structure
$J_1$, $\NN_1(X,Y)=0$ and
\begin{equation}
\label{eq:d2id2}
[J_1X, J_1Y]=0, \quad [J_1X,Y]=[J_1Y,X].
\end{equation}
Because of the integrability of the product  structure $J_2$,
$\NN_2(X,J_1X)=0$ and
\begin{equation}
\label{eq:d2id1} [X, J_1X]=[Y, J_1Y].
\end{equation}
From the Jacobi identity we have
$$
   \a(X,J_1X)=\b(X,J_1Y), \quad \a(J_1Y,X)=-\b(X,J_1X),
$$
and the bracket in $\g$ is determined by $c=\a(X,J_1X)$ and
$d=\b(X,J_1X)$ as follows:
$$
[X, J_1X]=[Y, J_1Y]=cX+dY,  \quad
[X, J_1Y]=[Y, J_1X]=dX+cY.
$$
Since $\dim \g'=2$,  $c^2-d^2\neq 0$ and we may choose
\begin{eqnarray*}
\xt = (c^2-d^2)^{-1}(cX+dY),\phantom{J_1 J_1} \qquad  &~&
\yt=(c^2-d^2)^{-
1}(dX+cY),\\
\zt = (c^2-d^2)^{-1}(cJ_1X-dJ_1Y), \qquad &~& \wt=
(c^2-d^2)^{-1}(-dJ_1X+cJ_1Y),
\end{eqnarray*}
and hence
\begin{eqnarray*}
~[\xt,\zt]=\xt, \qquad &~& [\xt,\wt]=\yt, \\
~[\yt,\zt]=\yt, \qquad &~& [\yt,\wt]=\xt,
\end{eqnarray*}
and we get algebra PHC7 for $a=1, b=0.$

{\bf Case iii):} $\g '$ is a totally null plane.  According to
Lemma \ref{le:2d_plane}  iii) we have to consider two
geometrically different cases.

In the first case we  can assume that   $J_2 |_{\g '} = 1$ and $\g
= \g ' + J_1\g '$ holds. If $(X, Y)$ is a basis of $\g'$ we have
\begin{equation*}
J_2X = X, \enskip J_2Y = Y,\enskip J_2J_1X = -J_1X, \enskip
J_2J_1Y = -J_1Y.
\end{equation*}
 One easily checks that the integrability of the
complex structure $J_1$ is equivalent to the relations
$$
[J_1X, J_1Y] = 0, \quad [X, J_1Y] = [Y, J_1X].
$$
It is interesting that the product structure $J_2$ is
automatically integrable. Hence, the possible non-null commutators
are
\begin{eqnarray*}
T' = [X, J_1X] &=& aX + bY,\\
Y' = [Y, J_1Y] &=& cX + dY,\\
X' = [X, J_1Y] &=& eX + f Y.
\end{eqnarray*}
The Jacobi identity is equivalent to the equations
\begin{equation}
\label{eq:jacobi}
 (e-d)X' +fY'-cT' = 0, \quad (a-f)X' +bY'-eT' = 0,
\end{equation}
or equivalently
\begin{equation*}
e(e-d)+ c(f-a) = 0, \quad ef = bc,\quad af-f^2+bd-be=0.
\end{equation*}
If $X'$ is a zero vector then we get the algebra PHC8. Suppose
that $X'$ is a non-zero vector. If $Y'$ or $T'$ is a zero vector
then we get an algebra PHC7 for $a=0=b.$ Suppose that none of the
vectors $X', Y', Z'$ is the zero vector. We can suppose that one
of the pairs $X',Y'$ and $X', T'$ is independent, say $X', T'.$ If
the vectors $X'$ and $Y'$ are collinear then we get the algebra
PHC7 for $a=0, b=1$. Finally, if the  both pairs $X', T'$ and $X',
Y'$ are independent then introduce a new basis $(X', Y', Z', W')$
satisfying
\begin{equation*}
Z' = \frac{1}{D}(fJ_1X - bJ_1Y), \quad W' = \frac{1}{D}(-eJ_1X+
aJ_1Y),
\end{equation*}
where $D = af - be\neq 0.$ In the new basis the commutator
relations take the very simple form
\begin{eqnarray*}
[X', Z'] &=& X', \enskip [X', W']=Y', \enskip [Y',Z']=Y',\\
 ~ [Y', W'] &=& \frac{fc-de}{D}X' + \frac{ad-bc}{D}Y'.
\end{eqnarray*}
Since $X'$ and $Y'$ are independent then $cf-de \neq 0$, that is,
$a \neq 0$ in the algebra PHC7.

In the second case we  can assume that   $(N_1, N_2)$ is a basis
of $\g '$ and $\g'$ is invariant with respect to $J_1,$ $J_2,$
$J_3$. Then a possible basis of $\g$ is
$$
N_1=J_1X+J_2X, \quad N_2=X-J_3X, \quad N_3=J_1X-J_2X, \quad
N_2=X+J_3X.
$$
We calculate the structures in terms of that basis:
$$
J_1N_1=-N_2, \quad J_1N_3=-N_4,
$$
$$
J_2N_1=N_2, \quad J_2N_3=-N_4,
$$
$$
J_3N_1=N_1, \quad J_3N_2=-N_2,\quad J_3N_3=-N_3, \quad J_3N_4=N_4.
$$
By the integrability of  $J_3$,
$$
J_3[N_1,N_4]=[N_1,N_4],\quad J_3[N_2,N_3]=-[N_2,N_3].
$$
Thus,
$$
[N_1,N_4]=\mu N_1,\quad [N_2,N_3]=\lambda N_2.
$$
The integrability of $J_1$ and $J_2$ is equivalent to
\begin{equation*}
0=-[N_2, N_4]-\lambda N_1 + \mu N_2 + [N_1, N_3]
\end{equation*}
After imposing the Jacobi identity this reduces to the algebra
PHC3.

{\bf Case iv):} the induced metric on $\g '$ is of rank $1$.
Denote by $N$ the  null vector belonging to $\g '$ (which is
unique up to a scaling constant).

According to  Lemma \ref{le:2d_plane}  iv) we can choose a product
structure $J_2$ such that for the basis $(X,N)$ of $\g '$ one has
\begin{equation}
N = J_1X - J_2X, \quad N \enskip \mbox{is null}.
\end{equation}
Then $(X,N, J_1X, J_1N)$ is a basis of $\g .$ One easily
calculates   the following relations
\begin{equation*}
J_2X = J_1X - N, \quad J_2N = J_1N.
\end{equation*}
The integrability of $J_1$ is equivalent to $\NN_1[X, N] = 0$,
i.e. to the relations
\begin{equation*}
[J_1X, J_1N]=0, \quad [X, J_1N] = [N, J_1X].
\end{equation*}
Since $(X, N, J_2X, J_2N)$ is a basis of $\g$ the integrability of
the product structure $J_2$ is equivalent to $\NN_2(X, N)=0$ which
gives  the condition
\begin{equation*}
[N, J_1N]=0.
\end{equation*}
The commutator relations now read
\begin{equation*}
[X, J_1X] = aX+bN,  \quad [X,J_1N]=cX+dN,
\end{equation*}
 where $a, b, c, d$ are unknown
coefficients. The Jacobi identity is now equivalent to the
following relations
\begin{equation*}
c = 0, \quad d(a-d) = 0.
\end{equation*}
The case $d = 0$ gives the algebra with $\dim \g ' =1$
 which we have already discussed. The
remaining case $a = d\neq 0,$ after the change
\begin{equation*}
\yt = N, \enskip \zt = J_1N, \enskip \xt = \frac{1}{a}X, \enskip
\wt = \frac{1}{a}J_1X - \frac{b}{a^2}J_1N,
\end{equation*}
takes the form
$$
  [\yt,\zt]=0,\enskip [\yt,\wt] =\yt, \enskip [\xt, \zt] = \yt,
  \enskip  [\xt,\wt] = \xt
$$
of the algebra PHC7 for $a=0=b$. \qed

\subsection{Case of solvable Lie algebra $\ggg$ with $\dim \ggg ' =
3$} \label{ssec:dim3}

\begin{Theorem}
\label{th:g3} Let  $\ggg$ be a $4$-dimensional solvable Lie
algebra  admitting a para-hypercomplex structure and with $\dim
\ggg '=3.$ If $\ggg$ has a nontrivial center it is algebra PHC6,
otherwise it is algebra PHC9 or PHC10.
\end{Theorem}
{\bf Proof: } If the algebra $\ggg$ is solvable then its derived
algebra $\ggg ' $ is nilpotent. Up to isomorphism the only
$3$-dimensional nilpotent Lie algebras are Abelian algebra and the
Heizenberg algebra generated by $X, Y$ and $Z$ with nonzero
commutator
$$
[X, Y] = Z.
$$
Let $\ggg$ be with trivial center, admitting  a para-hypecomplex
structure $(J_1, J_2)$ and let $\langle\cdot,\cdot\rangle$ be a
compatible inner product on $\ggg$. First, we discuss the case of
$\ggg ' $ being abelian.

 Suppose that $\ggg ' $ is nondegenerate subspace
and $X$ is normal vector of $\ggg '$. Then $|X|^2 \neq 0$ and $
\ggg ' = \RR \langle  J_1X, J_2X, J_3X \rangle .$ From the
integrability of $J_1$ and $J_2$ we have
$$
[X, J_\alpha J_\beta X]= J_\alpha[X,J_\beta X],
$$
for $\alpha, \beta\in{1,2,3},\ \alpha\neq\beta$.
Hence,
$[X, J_\alpha X]= \lambda J_\alpha$,
 and we get  the algebra PHC9 for $a = 0 = b.$
(the Lie algebra corresponding to the real hyperbolic spaces).
\vspace{0.5em}

Assume now  that $\ggg ' $ is degenerate subspace
 and $N$ is normal vector of $\ggg '.$ Then $|N|^2 = 0$ and $N\in \ggg '.$
According to  Lemma~\ref{le:2d_plane} iv) we can chose a
compatible structure $(J_1, J_2)$ such that $N = J_1X-J_2X$ for
any $X\in \ggg', \, |X|^2 \neq 0$. Since $J_1N$ is orthogonal to
$N$ we also have  $J_1N \in \ggg '$. Hence we may suppose that $
\ggg ' = \RR \langle  N, J_1N, X \rangle .$ Moreover the $(N,
J_1N, X, J_1X)$ is a basis of $\ggg.$ The integrability of $J_1$
and $J_2$ implies
$$
[J_1N,J_1X]=J_1[N,J_1X]=J_2[N,J_1X],
$$
 i.e. $ [N,J_1X]=d N$ and
 $[J_1N,J_1X]=d J_1N$, $d\neq 0$ what after scaling reduces to  algebra PHC9.

\vspace{0.5em} \noindent Now we turn to the case when  $\ggg '$ is
Heizenberg algebra. Let $\ggg ' = \RR \langle  X, Y, Z \rangle $
and $\ggg = \RR \langle   X, Y, Z , W \rangle$. One can easily
check that the center $\RR \langle Z \rangle$ is an ideal of
$\ggg$ and hence
$$
[W, Z ] = \lambda Z, \enskip \lambda \neq 0,
$$
no matter how the vector $W$ that does not belong to $\ggg '$ is
chosen. At the other side, independently of the choice of
non-central  vectors $X, Y \in \ggg '$ their commutator is always
in the  center, i.e.
$$
[X, Y] = \mu Z, \enskip \mu \neq 0.
$$
Here, $\mu \neq 0$ since $\ggg ' $ is not abelian and $\lambda
\neq 0$ since otherwise $Z$ would be a non-zero central element of
$\ggg .$ Hence, it remains to calculate the commutators $[W,X]$
and $[W, Y].$ This approach we use to prove the remaing part of
the theorem.

We consider the cases depending on degeneracy of $\ggg '$ with
respect to the induced compatible metric. Also there are different
subcases depending on the norm of a central element of $\ggg '.$

\vspace{0.5em} \noindent i) Suppose that $\ggg '$ is not
degenerated, and let $W$ be its normal vector. Denote by $Z = \xi
(\ggg ') $ a non-zero central element of $\ggg '.$ As an element
of $\ggg '$, $Z$ is orthogonal to $W.$ Now we have the following
cases.

{\bf $W$ and
 $Z$ have the same sign: } Using the Lemma
\ref{le:2d_plane} i) we may choose a compatible structure  $(J_1,
J_2)$ such that $Z = J_1W.$ Then the $(J_1W, J_2W, J_3W)$ is a
basis of $\ggg '$. After a simple calculation (and scaling) we get
the commutator relation:
$$
[W, J_1W] = 2 J_1W, \enskip [W, J_2W] =  J_2W, \enskip [W, J_3W] =
J_3W, \enskip [J_2W, J_3W ] = J_1W.
$$
That is a special form of algebra PHC10.

{\bf $W$ and $Z$ have  the opposite   sign: } Using  Lemma
\ref{le:2d_plane} ii) we may choose a compatible structure  $(J_1,
J_2)$ such that $Z = J_2W.$ Then the $(J_1W, J_2W, J_3W)$ is a
basis of $\ggg '$. After a simple calculation (and scaling) we get
the commutator relation:
$$
[W, J_1W] =  J_1W, \enskip [W, J_2W] = 2 J_2W, \enskip [W, J_3W] =
J_3W, \enskip [J_1W, J_3W ] = J_2W.
$$
That is again a special form of  algebra PHC10.

{\bf The center $Z$ of $\ggg '$ is a null vector: }  We have:
$|W|^2 \neq 0$, $|Z|^2 = 0$, $Z\perp X$, so using the Lemma
\ref{le:2d_plane} iv) we may choose a structure $(J_1, J_2)$ such
that
$$
N = Z = J_1W - J_2 W.
$$
Moreover there is a decomposition
$$
\ggg = \ggg ' \oplus \RR W = \RR \langle N, J_1W, J_3W \rangle
\oplus \RR W.
$$
Now we have
$$
[J_1W, J_3W] = \lambda N, \enskip [W, N] = \mu N, \enskip \lambda
, \mu \neq 0.
$$
After imposing the integrability condition for the structure
$(J_1, J_2)$ we get $\mu = 0$ what is a contradiction. Hence, this
case does not give a solution.

ii) Suppose that $\ggg '$ is degenerated, and let $N\in \ggg '$ be
its normal vector and $Z\in \ggg '$, a non-zero central element of
$\ggg '.$ We now discuss cases depending on the type of vector
$Z$.

{\bf $Z$ is a non null vector, $|Z|^2 \neq 0$:}
Let $X = Z.$
Consider the basis:
$$
\ggg = \RR \langle N, J_1N, X, J_1X \rangle, \enskip  \ggg '= \RR
\langle N, J_1N, X \rangle .
$$
Let $[N,J_1N]=\mu X$ and  $[J_1X,X]=\lambda X$. Then
$$
(J_1-J_2)[N,J_1X]=-\mu X.
$$
Thus, $\mu=0$ and $\g'$ is Abelian, what is again a contadiction.
\qed

{\bf $Z$ is a  null vector, $|Z|^2 = 0$:} According to Lemma
\ref{le:3d_plane}
  all null vectors of
$\ggg '$ are contained in two $2$-dimensional planes:
$$
\nnull (\ggg ') = \pi _1 \cup \pi_- = \RR \langle N, J_1N \rangle
\cup \{V | J_3 V = -V\}.
$$
We now study three possible cases $Z = N$, $Z \in \pi _-$ and $Z
\in \pi _1.$

{ ${\mathbf Z = N}$ {\bf (the normal to $\ggg '$ is a center of
$\ggg '$)}: } Then we have a decomposition:
$$
\ggg = \RR \langle N, J_1N, X, J_1X \rangle, \enskip \ggg ' = \RR
\langle N, J_1N, X \rangle.
$$
Because of the integrability of para-hypecomplex stucture
$(J_1,J_2)$ we have
$$
[J_1N, X] = \lambda N, \enskip [J_1X, N] = \mu N, \enskip
[J_1X,X]=aN+bJ_1N+cX \enskip \lambda , \mu \neq 0.
$$
The  Jacobi identity is equivalent to $c=\lambda$. After some
scaling we get the algebras PHC10.

  {\bf $Z \in\pi _-, Z \neq N$, ($Z$ is $-1$
eigenvector of $J_3$)}. Then $Z=aN+b(J_1N+2X)$ and we have  the
decomposition:
$$
\ggg = \RR \langle N, J_1N, Z, J_1Z \rangle, \enskip  \ggg '= \RR
\langle N, J_1N, Z \rangle.
$$
Due to the Heizenberg algebra structure of $\g'$ we  may assume
$$
[Z,J_1Z]=\mu Z, \quad [N,J_1N]=\lambda Z, \quad \mu,\lambda\neq 0.
$$
Because of the intergability of $J_1$ and $J_2$ we have
$$
[J_1N,J_1Z]=J_1[N,J_1Z]=J_2[N,J_1Z],
$$
and then
$$
[N,J_1Z]=\alpha N, \mbox{ and }  [J_1N,J_1Z]=\alpha J_1 N,
\alpha\neq 0.
$$
Now, by the Jacobi identity,
\begin{eqnarray*}
  [N, J_1Z]         & = & \alpha N,        \qquad     [Z, J_1 Z]=2\alpha Z,
\cr
  [J_1 N,J_1 Z]  & = & \alpha J_1 N,  \qquad     [Y, X] = \lambda Z,
\end{eqnarray*}
$\alpha, \lambda \neq 0$. After scaling it is a special case of
relations PHC10.

${\mathbf Z \in \mathbf \pi_1, Z = aN + J_1N, a\in \RR}.$ Consider
the decomposition
$$
\ggg = \RR \langle N, Z, X, J_1X \rangle, \enskip  \ggg '= \RR
\langle N, Z, X \rangle.
$$
Let $[N,X]=\mu Z$ and  $[J_1X,Z]=\lambda Z$. By the integrability,
$$
(J_1-J_2)[N,J_1X]=2\lambda Z-2\mu a N,
$$
what  implies $\lambda=0$, i.e. $Z$ is in the center of $\g$. That
is a contradiction.

\subsection{The proof of Theorem \ref{th:main}}
\label{ssec:exist}

According to the Levi decomposition theorem every Lie algebra
$\ggg$ decomposes into direct sum
$$
\ggg = \rrr \oplus \sss,
$$
where $\rrr$ is maximal solvable ideal (radical) and $\sss$ is
semisimple part. Since ${\mathfrak so}(3)$ and ${\mathfrak
sl}_2(\RR )$ are the only semisimple Lie algebras of dimension
less or equal to $4$, the only non-solvable Lie algebras of
dimension four are
$$
\RR \oplus {\mathfrak so}(3)\quad \mbox{and} \quad \RR \oplus
{\mathfrak sl}_2(\RR ).
$$
They both have a non-trivial center $\RR$, so from   Theorem
\ref{th:withcenter} we conclude that the unique non-solvable Lie
algebra  admitting a para-hypercomplex structure is $\RR \oplus
{\mathfrak sl}_2(\RR )$, i.e. PHC2. Solvable $4$-dimensional Lie
algebras  with nontrivial center and admitting a para-hypercomplex
structure are PHC1 and PHC3-PHC6 (Theorem \ref{th:withcenter}).
Solvable $4$-dimensional Lie algebras  with trivial center and
admitting a para-hypercomplex structure are PHC7-PHC10 (theorems
\ref{th:dim1}, \ref{th:dim2} and \ref{th:g3}).

It remains to prove that algebras PHC1-PHC10 posses  an integrable
para-hyper\-complex structure. We construct the structures below
and leave the reader to check   the integrability conditions for
$J_1$ and $J_2$ and the relations (\ref{eq:comm_rel}) by direct
calculation.

\noindent
 {\bf PHC1 and PHC2:}
$$
J_1Z = X, \enskip J_1Y = W, \quad J_2Z = Y, \enskip J_2X = -W.
$$

\noindent
 {\bf PHC3:}
\begin{eqnarray*}
J_1Z = X, \enskip J_1Y = W, &\\
J_2Z = W - Z, \enskip J_2X =X + Y,& \enskip J_2 Y = - Y, \enskip
J_2W = W.
\end{eqnarray*}

\noindent
 {\bf PHC4 and PHC5:}
\begin{eqnarray*}
J_1Z = W, \enskip J_1X = Y, &\\
J_2Z = W, \enskip J_2X =Y-Z,& \enskip J_2 Y = X + W.
\end{eqnarray*}

\noindent
 {\bf PHC7}
\begin{eqnarray*}
 J_1X = Z,  \enskip J_1Y = W,& \\
 J_2X = X, \enskip J_2Y = Y, & \enskip J_2Z = -Z, \enskip J_2W = -W
\end{eqnarray*}

\noindent
 {\bf PHC8:}
$$
J_1X = -Y, \enskip J_1Z = -W,\enskip  J_2X = Y, \enskip J_2Z = -W
$$

\noindent
 {\bf PHC6, PHC9 and PHC10:}
\begin{eqnarray*}
J_1Z = Y, \enskip J_1X = W, &\\
J_2Z = Y,  \enskip J_2 X = W - Z, \enskip & J_2W = X + Y . \qed
\end{eqnarray*}

\section{Comparisons with the work of Barberis}
\label{sec:comp}

In this section we compare our results with the classification of
hypercomplex structures in the paper of Barberis~\cite{B}.  We see
that there are many more $4$-dimensional Lie algebras with
para-hy\-per\-com\-plex structure  than Lie algebras with
hypercomplex structure.

Namely, we have the following.
\begin{Theorem}{(\cite{B})}
The only $4$-dimensional Lie algebras admitting an integrable
hypercomplex structure are:
\begin{itemize}
\item[(HC1)] $\g$ is abelian,
\item[(HC2)] $[X, Y]=W, \enskip [Y, W]=X, \enskip [W, X]=Y,$
\item[(HC3)] $[X,Z]=X, \enskip [X,W]=Y, \enskip [Y, Z]=Y, \enskip
[Y,W]=-Y,$
\item[(HC4)] $[W, X] = X, \enskip [W, Y] = Y, \enskip [W, Z] = Z$,
\item[(HC5)] $[W, X] = X,\enskip [W, Y] = \frac{1}{2}Y,\enskip [W, Z] =
\frac{1}{2}Z, \enskip [Z, Y] = X$.
\end{itemize}
\end{Theorem}

The Lie algebra HC2 is isomorphic to $\RR \oplus \frak{so}(3)$ and
it does not admit a para-hypercomplex structure. Its counterpart
admitting a para-hypercomplex (but not hypercomplex) structure is
algebra $\RR \oplus \frak{sl}(2)$ given by the relations PHC2.

No algebra $\ggg$ with $\dim \ggg ' = 1$ admits a hypercomplex
structure, while algebras PHC4 and PHC5 admit a para-hypercomplex
structure and satisfy $\dim \ggg ' = 1.$

The Lie algebra HC3 is isomorphic to $\frak{aff} (\CC )$ and it is
the only Lie algebra with $\dim \ggg ' = 2$ admitting a
hyper-complex structure. It also admits a para-hypercomplex
structure (PHC7 for $a = 1,$ $b = -1$).

The Lie algebra HC4 corresponds to real hyperbolic space $\RR
H^4$. It admits both hypercomplex and para-hypercomplex structure
(PHC9 for $a = 0 = b$).

Finally, the  Lie algebra HC5 corresponds to complex hyperbolic
space  $\CC H^2$. It admits both hypercomplex and
para-hypercomplex structure ( PHC10 for $c = 1$, $a=b=0$).




\begin{thebibliography}{10}

\bibitem{Sal1}
A.~Andrada, S~ Salamon, {\em Complex Product Structures on Lie
Algebras}, preprint math.DG/0305102.
\bibitem{B}
M.L.~Barberis, {\em  Hypercomplex Structures on Four-dimensional
Lie Groups,} Proc. of AMS, 128 (4) (1997), 1043--1054.
\bibitem{Ov}
G. Ovando, {\em Invariant complex structures on solvable real Lie
groups},
Manuscripta Math. 103(2000), 19--30.
\bibitem{Sal}
V. DeSmedt, S. Salamon,  {\em Anti-self-dual metrics on Lie
groups}, Proc. Conf. Integrable Systems and Differential Geometry,
Contemp. Math. 308(2002), 63--75.
\bibitem{S} Snow, J. E.,  {\em Invariant Complex Structures
on Four-dimensional Solvable Real Lie Groups}, Manuscripta Math.
66, (1990), 397--412.
\end{thebibliography}
\end{document}